\documentclass[aos,MSNbibl,citesort,dvips]{arximspdf}
\usepackage{graphicx}

\doi{10.1214/11-AOS946}
\volume{39}
\issue{6}
\pubyear{2011}
\firstpage{3369}
\lastpage{3391}

\makeatletter

\def\bsuffix #1{#1}

\newcommand{\con}{r}
\newcommand{\R}{\mathbb{R}}

\newtheorem{theorem}{Theorem}
\newtheorem{lemma}{Lemma}

\newproclaim{Remark}{Remark}

\makeatother

\begin{document}
\begin{frontmatter}

\title{Asymptotic optimality of the Westfall--Young permutation
procedure for multiple testing under~dependence\thanksref{T1}}
\runtitle{Optimality of the Westfall--Young permutation procedure}

\thankstext{T1}{N. Meinshausen and M. H. Maathuis contributed equally to this work.}

\begin{aug}
\author[A]{\fnms{Nicolai} \snm{Meinshausen}\corref{}\ead[label=nm]{meinshausen@stats.ox.ac.uk}},
\author[B]{\fnms{Marloes H.} \snm{Maathuis}\ead[label=mm]{maathuis@stat.math.ethz.ch}}\\
\and
\author[B]{\fnms{Peter} \snm{B\"uhlmann}\ead[label=pb]{buhlmann@stat.math.ethz.ch}}
\runauthor{N. Meinshausen, M. H. Maathuis and P. B\"uhlmann}
\affiliation{University of Oxford, ETH Z\"urich and ETH Z\"urich}
\address[A]{N. Meinshausen\\
Department of Statistics\\
University of Oxford\\
United Kingdom\\
\printead{nm}}
\address[B]{M. H. Maathuis\\
P. B\"uhlmann\\
Seminar f\"ur Statistik\\
ETH Z\"urich \\
Switzerland\\
\printead{mm}\\
\hphantom{E-mail: }\printead*{pb}} 
\end{aug}

\received{\smonth{6} \syear{2011}}
\revised{\smonth{11} \syear{2011}}

%
\begin{abstract}
Test statistics are often strongly dependent in large-scale multiple
testing applications. Most corrections for multiplicity are unduly
conservative for correlated test statistics, resulting in a loss of
power to detect true positives. We show that the Westfall--Young
permutation method has asymptotically optimal power for a broad class
of testing problems with a block-dependence and sparsity structure
among the tests, when the number of tests tends to infinity.
\end{abstract}

%

\end{frontmatter}

\section{Introduction}

We consider multiple hypothesis testing where the underlying tests are
dependent. Such testing problems arise in many applications, in
particular, in the fields of genomics and genome-wide association
studies \cite
{hirschhorn2005genome,mccarthy2008genome,dudoit2008multiple}, but also
in astronomy and other fields \cite
{liang02statistical,meinshausen04estimating}. Popular multiple-testing
procedures include the Bonferroni--Holm method \cite{holm79simple}
which strongly controls the family-wise error rate (FWER), and the
Benjamini--Yekutieli procedure \cite{benjamini01control} which controls
the false discovery rate (FDR), both under arbitrary dependence
structures between test statistics. If test statistics are strongly
dependent, these procedures have low power to detect true positives.
The reasons for this loss of power are well known: loosely speaking,
many strongly dependent test-statistics carry only the information
equivalent to fewer ``effective'' tests. Hence, instead of correcting
among many multiple tests, one would in principle only need to correct
for the smaller number of ``effective'' tests. Moreover, when
controlling some error measure of false positives, an oracle would only
need to adjust among the tests
corresponding to true negatives. 
In large-scale
sparse multiple testing situations, this latter issue is
usually less important since the number of true positives is typically
small, and the number of true negatives is close to the
overall number of tests.

The dependence among
tests can be taken into account by using the permuta\-tion-based
Westfall--Young method \cite{westfall93resampling}, already used widely
in practice (e.g., \mbox{\cite{cheung2005mapping,winkelmann2007genome}}).
Under the assumption of subset-pivotality
(see Section \ref{secsinglestepWY} for a definition),
this method strongly controls the FWER under any kind of dependence
structure \cite{WestfallYoung89}.

In this paper we show that the Westfall--Young
permutation method is an optimal procedure in the following sense. We introduce
a single-step oracle
multiple testing procedure, by defining a single threshold such
that all hypotheses with $p$-values below this threshold are rejected
(see Section \ref{secsinglestep}). The
oracle threshold is the largest threshold that still guarantees the desired
level of the testing procedure. The oracle threshold is unknown in practice
if the dependence among test statistics and the set of true null hypotheses
are
unknown.
We show that the single-step Westfall--Young threshold
approximates
the
oracle threshold for a broad class of testing problems with a~block-dependence
and sparsity structure among the tests, when the number of tests tends
to infinity.
Our notion of asymptotic optimality relative to an oracle threshold is
on a
general level and for any specified test statistic. The power of a multiple
testing procedure depends also on the data generating distribution and the
chosen individual test(s): we do not discuss this aspect here. Instead, our
goal is to analyze optimality once the individual tests have been specified.

Our
optimality result has an immediate consequence for large-scale
multiple testing: it is not possible to improve on the power of the
Westfall--Young permutation method while still controlling the FWER
when considering
single-step multiple testing procedures for
a large number of tests and assuming only a
block-dependence and sparsity structure
among
the tests (and no additional modeling
assumptions about the dependence or clustering/grouping).
Hence, in such situations, there is no need to consider ad-hoc proposals
that are sometimes
used
in practice, at least when taking the viewpoint
that multiple testing adjusted $p$-values
should be as model free as possible.

\subsection{Related work}
There is a small but growing literature on optimality in multiple
testing under dependence. Sun and Cai \cite{sun2009large} studied and proposed
optimal decision procedures in a two-state hidden Markov model, while
Genovese et al. \cite{genovese2006false} and Roeder and Wasserman
\cite{roeder2009genome} looked at the
intriguing possibility of incorporating prior information by $p$-value
weighting. The effect of correlation between test statistics on the
level of FDR control was studied in Benjamini and Yekutieli \cite
{benjamini01control} and
Benjamini et al. \cite{benjamini2006adaptive}; see also
Blanchard and Roquain \cite{blanchard2009adaptive}
for FDR control under dependence. Furthermore,
Clarke and Hall \cite{clarke2009robustness} discuss the effect of
dependence and
clustering when using ``wrong'' methods based on independence
assumptions for controlling the (generalized) FWER and FDR. The effect
of dependence on the power of Higher Criticism was examined in
Hall and Jin \cite{hall2008properties,hall2010innovated}. Another
viewpoint is
given by Efron \cite{efron2007correlation}, who proposed a novel empirical
choice of an appropriate null distribution for large-scale significance
testing. We do not propose new methodology in this manuscript but study
instead the asymptotic optimality of the widely used Westfall--Young
permutation method \cite{westfall93resampling} for dependent test
statistics.

\section{Single-step oracle procedure and the Westfall--Young method}

After introducing some notation, we define our notion of
a single-step oracle
threshold and describe the Westfall--Young permutation method.

\subsection{Preliminaries and notation}\label{secnotation}
Let $W$ be a data matrix containing $n$ independent realizations of an
$m$-dimensional
random variable $X = (X_1,\ldots,\allowbreak X_m)$ with distribution $P_m$ and possibly
some additional deterministic response variables $y$.\vspace*{12pt}

\textit{Prototype of data matrix $W$.}
To make this more concrete, consider the following setting that fits the
examples described in Section \ref{secexamples}. Let $y$ be a
deterministic variable, and allow the distribution of $X=X_y$ to depend on
$y$. For each value $y^{(i)}$, $i=1,\ldots,n$,
we observe an independent sample $X^{(i)}=(X^{(i)}_1,\ldots,X^{(i)}_m)$
of $X=X_{y^{(i)}}$. We then define $W$ to be an
$(m+1)\times n$-dimensional
matrix by setting $W_{1,i}=y^{(i)}$ for $i=1,\ldots,n$ and
$W_{j+1,i} =
X_j^{(i)}$ for $j=1,\ldots,m$ and $i=1,\ldots,n$. Thus, the first row of
$W$ contains
the $y$-variables, and the $i$th column of $W$ corresponds to the
$i$th data sample~$(y^{(i)},X^{(i)})$.\vspace*{12pt}

Based on $W$, we want to test $m$ null hypotheses $H_j$,
$j=1,\ldots,m$, concerning
the $m$ components $X_1,\ldots,X_m$ of $X$. For concrete examples, see
Section \ref{secexamples}.
Let $I(P_m) \subseteq\{1,\ldots,m\}$ be the indices of the true null
hypotheses,
and let $I'(P_m)$
be the indices of the true alternative
hypotheses, that is, $I'(P_m) = \{1,\ldots,m\} \setminus I(P_m)$.
Let $P_0$ be
a distribution under the complete null hypothesis,
that is, $I(P_0)=\{1,\ldots,m\}$. We denote the class of all
distributions under the complete null hypothesis by $\mathcal P_0$.

Suppose that the same test is applied for all hypotheses, and let
$S_n\subseteq[0,1]$ be the set of
possible $p$-values this test can take. Thus, $S_n=[0,1]$ for $t$-tests and
related approaches, while
$S_n$ is discrete for permutation tests and rank-based tests. Let
$p_j(W)$, $j=1,\ldots,m$, be the $p$-values for the $m$ hypotheses, based
on the chosen test and the data $W$.

\subsection{Single-step oracle multiple testing procedure}
\label{secsinglestep}

Suppose that we knew the true set of null hypotheses $I(P_m)$ and the
distribution of $\min_{j\in I(P_m)} p_{j}(W)$ under $P_m$ (which is of
course not true in practice).
Then we could define the following single-step
oracle multiple testing procedure:
reject $H_j$ if
$p_j(W) \le c_{m,n}(\alpha)$, where $c_{m,n}(\alpha)$ is the
$\alpha$-quantile of $\min_{j\in I(P_m)}p_j(W)$ under $P_m$.
%
\begin{equation}\label{eqoracle}
c_{m,n}(\alpha) = \max\Bigl\{ s \in S_n \dvtx P_m\Bigl( \min_{j\in I(P_m)}
p_j(W)\le s\Bigr) \le\alpha\Bigr\}.
\end{equation}
Throughout, we define the maximum of the empty set to be zero,
corresponding to a threshold
$c_{m,n}(\alpha)$
that leads to zero rejections.

This oracle
procedure
controls the FWER at level $\alpha$, since, by definition,
\begin{eqnarray*}
&& P_m \bigl( H_j \mbox{ is rejected for at least one } j\in I(P_m) \bigr)
\\
&&\qquad = P_m\Bigl( \min_{j\in I(P_m)} p_{j}(W) \le c_{m,n}(\alpha)\Bigr) \le
\alpha,
\end{eqnarray*}
and it is optimal in the sense that values $c \in S_n$ with $c >
c_{m,n}(\alpha)$ no longer control the FWER at level $\alpha$.

\subsection{Single-step Westfall--Young multiple testing
procedure}\label{secsinglestepWY}

The Westfall--Young permutation method
is based on the idea that under the complete null hypothesis, the
distribution of $W$ is invariant under a certain group of
transformations~$\mathcal G$,
that is, for every $g\in\mathcal G$, $gW$ and $W$ have the same distribution
under $P_{0}\in\mathcal P_0$.
Romano and Wolf \cite{RomanoWolf05} refer to this as the
``randomization hypothesis.'' In
the sequel, $\mathcal G$ is the collection of all permutations~$g$ of
$\{1,\ldots,n\}$, so that the number of elements $|\mathcal G|$
equals $n!$.\vspace*{12pt}

\textit{Prototype permutation group $\mathcal G$ acting on the
prototype data matrix $W$.} In the examples in Section
\ref{secexamples}, $W$ is a prototype data matrix as described in
Section \ref{secnotation}. The prototype permutation $g \in\mathcal G$
leads to a matrix $gW$ obtained by permuting the \textit{first} row of
$W$ (i.e., permuting the $y$-variables). For all examples in
Section~\ref{secexamples}, under the complete null hypothesis $P_0
\in\mathcal P_0$, the distribution of $gW$ is then identical to the
distribution of $W$ for all $g\in\mathcal G$, so that the randomization
hypothesis is satisfied. We suppress the dependence of $|\mathcal{G}|$
on the sample size $n$ for notational simplicity.\vspace*{12pt}

%

The single-step Westfall--Young critical value is a random variable,
defined as follows:
\begin{eqnarray*}
\hat c_{m,n}(\alpha) & = & \max\biggl\{ s \in S_n \dvtx\frac{1}{|\mathcal{G}|}
\sum_{g\in\mathcal G} 1\Bigl\{ \min_{j=1,\ldots,m} p_{j}(gW) \le s \Bigr\}
\le\alpha\biggr\} \\
& = &\max\Bigl\{ s \in S_n \dvtx P^*\Bigl(\min_{j=1,\ldots,m} p_{j}(W) \le s\Bigr)
\le\alpha\Bigr\},
\end{eqnarray*}
where $1\{\cdot\}$ denotes the indicator function,
and $P^*$ represents the permutation distribution
%
\begin{equation} \label{eqPstar}
P^*\bigl( f(W) \le x\bigr) = \frac{1}{|\mathcal{G}|} \sum_{g\in\mathcal G} 1\{
f(gW) \le x \}
\end{equation}
for any function $f(\cdot)$ mapping $W$ into $\mathbb{R}$.
In other words,
$\hat c_{m,n}(\alpha)$ is the $\alpha$-quantile of the permutation
distribution of $\min_{j=1,\ldots,m}p_j(W)$.
Our main result (Theorem~\ref{thoptimal}) shows that under some
conditions, the Westfall--Young threshold $\hat{c}_{m,n}(\alpha)$
approaches
the oracle threshold $c_{m,n}(\alpha)$.

It is easy to see that the Westfall--Young permutation method provides weak
control of the
FWER, that is, control of the FWER
under the complete null hypothesis. Under the assumption of
subset-pivotality, it also provides strong control of the FWER
\cite{westfall93resampling}, that is, control of the FWER under any set
$I(P_m)$ of true null hypotheses.
Subset-pivotality
means that the distribution of $\{p_j(W) \dvtx j\in K\}$ is identical under
the restrictions $\bigcap_{j\in K} H_j$ and $\bigcap_{j\in I(P_0)} H_j$ for all
possible subsets $K \subseteq I(P_m)$ of true null
hypotheses. Subset-pivotality is not a necessary condition for strong
control; see, for example,
Romano and Wolf \cite{RomanoWolf05}, Westfall and Troendle
\cite{westfall2008multiple} and
Goeman and Solari \cite{GoemanSolaris10}.

\section{Asymptotic optimality of Westfall--Young}\label{secoptimality}

We consider the framework where the number of hypotheses $m$ tends to
infinity. This framework is suitable for high-dimensional settings
arising, for example, in
microarray experiments or genome-wide association studies.

\subsection{Assumptions}\label{secassumptions}

(A1) Block-independence: the $p$-values of all true null hypotheses
adhere to a block-independence structure 
that is preserved under permutations in $\mathcal G$.
Specifically, there exists a partition
$A_1,\ldots,A_{B_m}$ of $\{1,\ldots,m\}$ such that 
for any pair of permutations $g,g' \in\mathcal G$,
\[
\min_{\tilde{g}\in\{g,g'\}} \min_{j\in A_b\cap I(P_m)} p_j(\tilde{g}W)
,\qquad b=1,\ldots,B_m,
\]
are mutually independent
under $P_m$.
Here, the number of blocks is denoted by $B=B_m$.
[We assume without loss of generality
that $A_b\cap I(P_m)\neq\varnothing$ for all $b=1,\ldots,B$,
meaning that there is at least one true null hypothesis in each block;
otherwise, the condition would be
required
only for blocks with $A_b \cap
I(P_m) \neq\varnothing$.]\vspace*{-5pt}

\begin{longlist}[{(A2)}]
\item[{(A2)}] Sparsity: the number of alternative hypotheses that are
true under~$P_m$ is small compared to the number of blocks, that is,
$|I'(P_m)| / B_m \to0$ as $m\to\infty$.
\item[{(A3)}] Block-size: the 
maximum size of a block, $m_{B_m}:= {\max_{b=1,\ldots,B_m}} |A_b|$, is of
smaller order than the square root of the number of blocks, that is,
$m_{B_m} =o(\sqrt{B_m})$ as $m\to\infty$.
\item[(B1)] Let
$G$ be a random permutation taken uniformly from $\mathcal{G}$.
Under~$P_m$,
the joint
distribution of $\{p_j(W) \dvtx j\in I(P_m)\}$
is identical to the joint
distribution of $\{p_j(GW) \dvtx j\in I(P_m)\}$.
\item[{(B2)}] Let $P^*$ be the permutation distribution in
(\ref{eqPstar}). There exists a constant $\con<\infty$ such that for
$s = c_{m,n}(\alpha) \in S_n$ and all $W$,
%
\begin{equation}\label{eqT1} \con^{-1} s \le P^*\bigl( p_j(W)
\le s\bigr) \le\con s \qquad\mbox{for all } j=1,\ldots,m.
\end{equation}
\item[{(B3)}] The $p$-values corresponding to true null hypotheses are
uniformly distributed; that is, for all $j\in I(P_m)$ and $s\in S_n$,
we have $P_m(p_j(W) \le s) = s$.
\end{longlist}

A sufficient condition for the block-independence assumption (A1) is
that for every fixed pair of permutations
$g, g' \in\mathcal G$ the blocks of random variables
$\{p_j(gW), p_j(g'W)\dvtx j\in A_b \cap I(P_m)\}$ are mutually independent for
$b=1,\ldots,B_m$. This condition is implied by block-independence of
the $m$
last rows of the prototype $W$ for
the examples discussed in Section \ref{secexamples} and for the prototype
$\mathcal G$ as in Section \ref{secsinglestepWY}. The block-independence
assumption
captures an essential characteristic
of large-scale testing problems: a test statistic is often
strongly correlated with a
number of other test statistics but
not at all with the remaining tests.

The sparsity assumption (A2)
is appropriate in many contexts. Most
genome-wide association studies, for example,
aim to discover just a few locations on the genome that are associated with
prevalence of a certain disease \cite
{kruglyak1999prospects,marchini2005genome}.
Furthermore, assumption (A3) requiring that the range of (block-)
dependence is not too
large, which seems reasonable in genomic applications: for example,
when having many different groups of genes (e.g., pathways), each of
them not too large in cardinality, a block-dependence structure seems
appropriate.

We now consider assumptions (B1)--(B3), supposing that we work with a~prototype data matrix $W$ and a prototype permutation group $\mathcal
G$ as described in Sections \ref{secnotation} and
\ref{secsinglestepWY}. Assumption (B1) is satisfied if each $p$-va\-lue~$p_j(W)$ only depends on the $1$st and $(j+1)$th rows of $W$. Moreover,
subset-pivotality is satisfied in this setting. Assumption (B3) is
satisfied for any test with valid type I error control. Assumption (B2)
is fulfilled with $r=1$ if for all $W$
%
\begin{equation}\label{eqcondprob}
P_G\bigl( p_j (GW) \le s |W \bigr) =s,\qquad
j=1,\ldots,m, s\in S_n,
\end{equation}
where $P_G$ is the probability with respect to a random permutation $G$
taken uniformly from $\mathcal{G}$, so that the left-hand side of
(\ref{eqcondprob})
equals
$ P^*( p_j(W) \le s)$ in~(\ref{eqT1}). Note that assumptions
(B1) and (B3) together imply that
%
\begin{equation}\label{eqtest}
P_{m,G}\bigl( p_j (GW) \le s\bigr) =s,\qquad j\in
I(P_m), s\in S_n,
\end{equation}
where the probability $P_{m,G}$ is with respect to a random draw of the
data~$W$,
and a random permutation $G$ taken uniformly from $\mathcal{G}$. Thus,
assumption~(B2) holds if
(\ref{eqtest}) is true for all $j=1,\ldots,m$ when conditioned on the
observed data.
Section \ref{secexamples} discusses
three concrete examples that satisfy assumptions~\mbox{(B1)--(B3)} and
subset-pivotality.
\begin{Remark*}
For our theorems in Section \ref{secmainresult}, it would be
sufficient if
(\ref{eqT1}) were holding only with probability converging to 1 when
sampling a random~$W$, but we leave a deterministic bound since it is
easier notationally, the extension is direct and we are mostly interested
in rank-based and conditional tests for which the deterministic bound
holds.
\end{Remark*}

\subsection{Examples}\label{secexamples}

We now give three examples that satisfy assump-\break tions~\mbox{(B1)--(B3)}, as well as
subset-pivotality.
As in Section \ref{secnotation}, let $y$ be a~deterministic scalar
class variable and $X=(X_1,\ldots,X_m)$ an $m$-dimensional vector of
random variables, where the distribution of $X=X_y$ can depend on~$y$. Let
the prototype data matrix $W$ and the prototype group of permutations
$\mathcal G$ be defined
as in Sections \ref{secnotation} and \ref{secsinglestepWY},
respectively.
In all examples, we work with tests with valid type I error control,
and each $p$-value $p_j(W)$ only depends on the
$1$st and $(j+1)$th rows of $W$. Hence, assumptions~(B1),~(B3) and~%
subset-pivotality are satisfied, and we focus on
assumption (B2) in the remainder.

%
%

For the examples in Sections \ref{seclocationshift} and \ref
{secmarginal}, we assume that there exists a $\mu(y)\in\mathbb{R}^m$
and an $m$-dimensional random variable $Z=(Z_1,\ldots,Z_m)$ such that
%
\begin{equation}\label{eqXy}
X=X_y = \mu(y) + Z .
\end{equation}
We omit the dependence of $X=X_y$ on $y$ in the following for
notational simplicity.

\subsubsection{Location-shift models}\label{seclocationshift}
We consider two-sample testing problems for location shifts, similar to
Example 5 of Romano and Wolf \cite{RomanoWolf05}.
Using the notation in (\ref{eqXy}), $y\in\{1,2\}$ is a binary class
variable, and the marginal distributions of $Z$ are assumed to have a median
of zero.

We are
interested in testing the
null hypotheses
\[
H_{j}\dvtx \mu_j(1) = \mu_j(2),\qquad j=1,\ldots,m,
\]
versus the corresponding two-sided
alternatives,
\[
H'_{j}\dvtx \mu_j(1) \neq\mu_j(2),\qquad j=1,\ldots,m.
\]
%

We now discuss location-shift tests that satisfy assumption (B2).
First, note that all permutation tests satisfy (B2) with $r=1$, since the
$p$-values in a permutation test are defined to fulfill $P^*( p_j(W)
\le s)
=s$ for all $s\in S_n$. Permutation tests are often recommended in
biomedical research \cite{ludbrook1998permutation} and other large scale
location-shift testing applications due to their robustness with
respect to the underlying
distributions. For example, one can use the Wilcoxon test. Another example
is a ``permutation $t$-test'': choose the $p$-value $p_j(W)$ as the
proportion of
permutations for which the absolute value of the $t$-test statistic is
larger than or equal to the
observed absolute value of the $t$-test statistic for $H_j$.
Then condition (B2) is
fulfilled with $r=1$ with the added advantage that inference is exact, and
the type I error is guaranteed even if the distributional Gaussian
assumption for the $t$-test is not fulfilled
\cite{good2000permutation}. Computationally,
such a ``permutation $t$-test'' procedure seems to involve two rounds of
permutations: one for the computation of the marginal $p$-value and one for
the Westfall--Young method; see (\ref{eqPstar}). However,
the computation of the marginal permutation $p$-value can be inferred
from the
permutations in the Westfall--Young method, as in
Meinshausen \cite{meinshausen03false}, and just a single round of
permutations is thus
necessary.

\subsubsection{Marginal association}\label{secmarginal}
Suppose that we have a continuous variable~$y$ in
formula (\ref{eqXy}). Based on
the observed data, we want to test the null hypotheses of no association
between variable $X_j$ and $y$, that is,
\[
H_{j}\dvtx \mu_j(y) \mbox{ is constant in } y,\qquad j=1,\ldots,m,
\]
versus the corresponding two-sided alternatives. A special case is the test
for linear marginal association, where the functions $\mu_j(y)$ for
$j=1,\ldots,m$ are assumed to be of the form $\mu_j(y) = \gamma_j +
\beta_j
y$, and the test of no linear marginal association is based on the null
hypotheses
\[
H_{j}\dvtx \beta_j = 0,\qquad j=1,\ldots,m.
\]

Rank-based correlation test like Spearman's or Kendall's correlation
coefficient are examples of tests that fulfill assumption (B2).
Alternatively, a~``permutation
correlation-test'' could be used, analogous to the ``permutation
$t$-test'' described
in Section \ref{seclocationshift}.

\subsubsection{Contingency tables}\label{seccontingency}
Contingency tables are our final example. Let $y\in\{1,2,\ldots,K_y\}$
be a class variable with $K_y$ distinct values.
Likewise, assume that the random variable $X$ is discrete and that each
component of $X$ can take $K_x$ distinct values, $X=(X_1,\ldots
,X_m)\in
\{1,2,\ldots,K_x\}^m$.

As an example, in many genome-wide association studies, the variables
of interest are
single nucleotide polymorphisms (SNPs). Each SNP $j$ (denoted by $X_j$)
can take three distinct values, in general, and it is of interest to see
whether there is a relation between the occurrence rate of these categories
and a category of a
person's
health status $y$
\cite{kruglyak1999prospects,goode2002effect,bond2005single}.

Based on the observed data, we want to test the null hypothesis for
$j=1,\ldots,m$ that
the distribution of $X_j$ does not depend on $y$,
\[
H_{j}\dvtx P(X_j=k|y) = P(X_j=k) \qquad\mbox{for all } k\in
\{1,\ldots,K_x\}\mbox{ and } y\in\{1,\ldots,K_y\}.
\]
The available data for hypothesis $H_j$ is contained in the $1$st and
\mbox{$(j+1)$th} rows of~$W$. These data can be summarized in a contingency
table and Fisher's exact test can be used.
Since the test is conditional on the marginal distributions, we have
that $
P( p_j (GW) \le s |W ) =s$ for a random permutation $G\in\mathcal G$ and
(B2) is fulfilled with $r=1$.

\subsection{Main result}\label{secmainresult}

We now look at the properties of the Westfall--Young permutation method and
show asymptotic optimality in the sense that, with probability
converging to 1 as the
number of tests increases, the estimated
Westfall--Young threshold $\hat c_{m,n}(\alpha)$
is at
least as large as the optimal oracle threshold $c_{m,n}(\alpha-\delta)$,
where $\delta>0$ can be arbitrarily small.
This implies that the power of the Westfall--Young permutation method
approaches the power of the oracle test, while providing strong control of
the FWER under subset-pivotality \cite{westfall93resampling}.
All proofs are given in Section \ref{secproofs}.\vspace*{-2pt}
%
\begin{theorem}\label{thoptimal}
Assume
\textup{(A1)--(A3)} and \textup{(B1)--(B3)}.
Then for any $\alpha\in(0,1)$ and any $\delta\in(0,\alpha)$
%
\begin{equation}\label{eqthoptimal}
P_m\{ \hat c_{m,n}(\alpha) \ge c_{m,n}(\alpha-\delta) \} \to1
\qquad\mbox{as } m\to\infty.
\end{equation}
\end{theorem}

We note that
the sample size $n$
can be fixed and does not need to
tend to
infinity. However, if the range of $p$-values $S_n$ is discrete, the sample
size must increase with $m$ to avoid a trivial result
where the oracle
threshold \mbox{$c_{m,n}(\alpha-\delta)$}
vanishes; see also Theorem \ref{thoptimaldiscrete} where this
is made explicit for the Wilcoxon test in the location-shift model of
Section \ref{seclocationshift}.

Theorem \ref{thoptimal} implies that the actual level of the
Westfall--Young procedure converges
to the desired level (up to possible discretization effects; see
Section \ref{sectiondiscrete}).
To appreciate the statement in Theorem \ref{thoptimal} in terms of power
gain, consider a simple example. Assume that the $m$ hypotheses form
$B$ blocks. In the most
extreme scenario, test statistics are perfectly dependent
within each block. In such a scenario, the
oracle threshold
(\ref{eqoracle}) for each individual $p$-value is then
\[
1- \sqrt[B]{1-\alpha} ,
\]
which is larger than, but very closely approximated by $\alpha/B$ for
large values of~$B$. Thus, when controlling the FWER at level $\alpha$,
hypotheses can be rejected when their $p$-values are less than $1-
\sqrt[B]{1-\alpha}$ and certainly when their $p$-values are less than
$\alpha/B$. However, the value of $B$ and the block-dependence
structure between hypotheses are unknown in practice.
With a~Bonferroni correction for
the FWER at level $\alpha$,
hypotheses can be rejected when their $p$-values are less than $\alpha/m$.
If $m\gg B$,
the power loss compared to the
procedure with the oracle threshold is substantial, since the Bonferroni
method is really controlling at an effective level of size $\alpha B/m$
instead of $\alpha$.
Theorem \ref{thoptimal}, in contrast, implies that the effective level
under the Westfall--Young procedure converges to the desired level (up
to possible discretization effects).\vadjust{\goodbreak}

\subsection{Discretization effects with Wilcoxon test}
\label{sectiondiscrete}
We showed in the last section that the Westfall--Young permutation
method is
asymptotically equivalent to the oracle threshold under the made
assumptions.
In this section we look in more detail at the difference between the nominal
and effective levels
of the oracle multiple testing procedure.
Controlling at nominal level~$\alpha$, the effective 
oracle level is defined as
%
\begin{equation}\label{eqalpha-}
\alpha_- = P_m \Bigl\{ \min_{j\in I(P_m) } p_j(W) \le c_{m,n}(\alpha)
\Bigr\}.
\end{equation}
By definition, $\alpha_-$ is less than or equal to $\alpha$.
We now examine under which assumptions the effective level
$\alpha_-$ can be replaced by the nominal level $\alpha$.
As a concrete example, we work with the following 
assumptions:
\begin{longlist}[{(A3$'$)}]
\item[{(W)}] The test is a two-sample Wilcoxon test with equal sample
sizes $n_1=n_2=n/2$, applied to a location-shift model as defined in
Section \ref{seclocationshift}.\vspace*{1pt}

\item[{(A3$'$)}] Block-size:
the maximum size of a block satisfies $m_B =O(1)$ as \mbox{$m\to\infty$}.
\end{longlist}
The restriction to equal sample sizes in (W) is only for technical
simplicity.
We then obtain the following result about the discretization
error.\vspace*{-2pt}
%
\begin{theorem}\label{thoptimaldiscrete}
Assume \textup{(W)}. Then the oracle critical value $c_{m,n}(\alpha)$ is
strictly positive when
\[
n \ge2 \log_2(m/\alpha) + 2.
\]
%
When assuming in addition \textup{(A1), (A2)} and \textup{(A3$'$)}, then
the results of Theorem \ref{thoptimal} hold, and, for any $\alpha\in
(0,1)$, we have
\[
\alpha_- \to\alpha
\]
as $m,n\to\infty$ such that $n/\log(m) \to\infty$.\vspace*{-2pt}
\end{theorem}

The first result in Theorem \ref{thoptimaldiscrete} says that the
oracle critical value for the test defined in (W) is nontrivial, even
when the number of tests grows almost exponentially with sample size.
Hence, in this setting the result from Theorem \ref{thoptimal} still
applies in a nontrivial way.

The second result in Theorem \ref{thoptimaldiscrete} gives sufficient
criteria for the effective oracle level $\alpha_-$ to converge to
$\alpha$. It is conceivable that this result can also be obtained under
a milder assumption than (A3$'$), but this requires a detailed study of
the Wilcoxon $p$-values, and we leave this for future work. The main
takeaway message is that discreteness of the $p$-values does not change
the optimality result fundamentally.

\section{Empirical results}\label{secempirical}

The power of the Westfall--Young procedure has already been examined
empirically in several studies. Westfall et al. \cite{wezay02}
includes
a comparison with the Bonferroni--Holm method, reporting a gain in power
when using the Westfall--Young procedure. Its focus is on ``genetic
effects in association\vadjust{\goodbreak}
studies,'' including genotype (SNP-type) analysis and also gene
expression microarray
analysis.
Becker and Knapp \cite{bekna04} apply the Westfall--Young permutation
procedure and report substantial gain in power over Bonferroni
correction in the context of haplotype analysis.
Yekutieli and Benjamini \cite{yebe99} and Reiner et al. \cite{reyebe03}
discuss the gain of resampling in
terms of power, although their focus is mainly on FDR controlling methods.
Also Dudoit et al. \cite{dudoit2003mht} and Ge et al.
\cite{GeEtAl03} report that resampling-based methods such as the
Westfall--Young permutation procedure have clear
advantages in terms of power.

Here, we look at a few simulated examples to study the finite-sample
properties and compare with the asymptotic results of Theorem
\ref{thoptimal}. Data for a two-sample location-shift model as in
Section \ref{seclocationshift} with $m$ hypotheses and equal sample
sizes of $50$ are generated from a multivariate Gaussian distribution
with unit variances and (i) a Toeplitz correlation matrix with
correlations $\rho_{i,j}= \rho^{|i-j|}$ for some $\rho\in
\{0.95,0.975,0.99\}$ and $1\le i,j\le m$ and (ii)~a~block covariance
model, where correlations within all blocks (of size~50 each) are set
to the same $\rho\in\{0.6,0.75,0.9\}$ and to 0 outside of each block.
Ten alternative hypotheses are picked at random from the first 100
components by applying a shift of 0.75 whereas all remaining $m-10$
null-hypotheses correspond to no shift.

\begin{figure}

\includegraphics{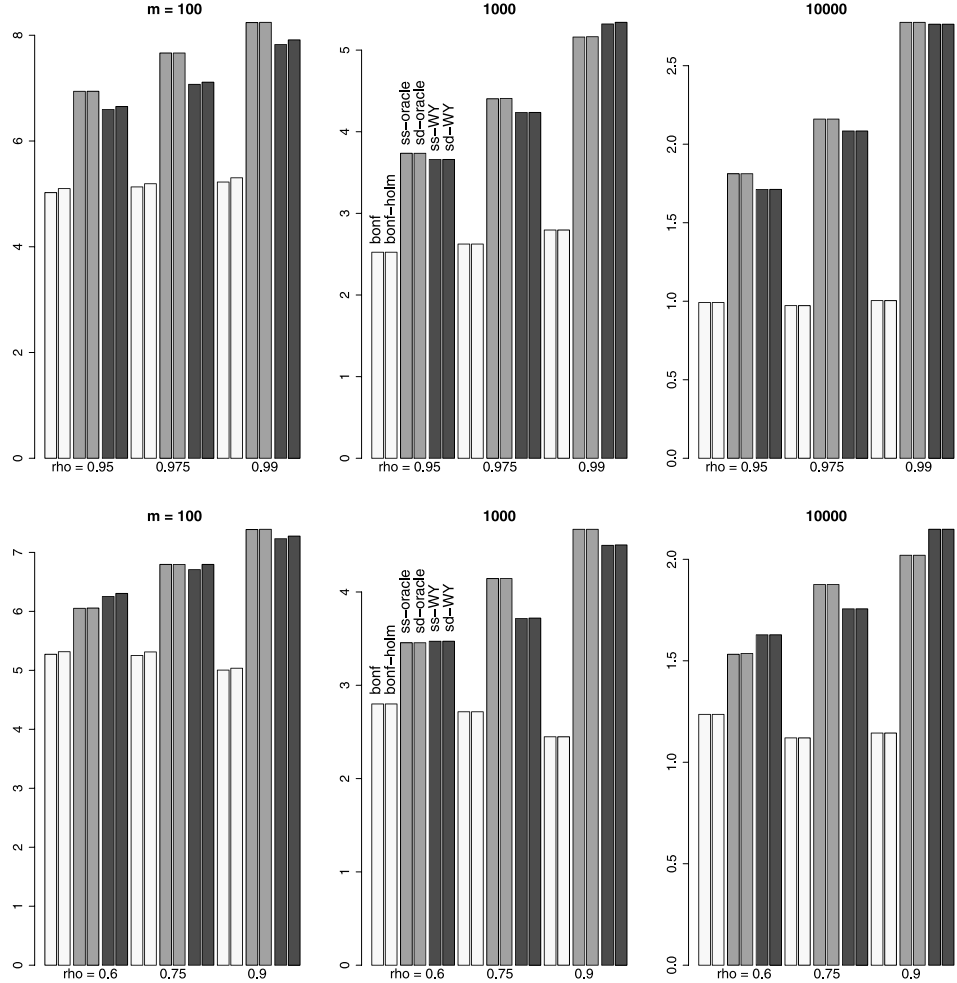}

\caption{Average number of true positives for the
Toeplitz model (top row) and the block model (bottom row). The number
of hypotheses is increasing from $m=100$ (left panel) to $m=10\mbox{,}000$
(right panel) and the correlation parameter is varied as $\rho
=0.95,0.975, 0.99$ in the Toeplitz model and $\rho=0.6,0.75, 0.9$ in
the block model. The results are shown for the Bonferroni correction
(single-step and corresponding step-down Bonferroni--Holm; light
color); the oracle procedure (single-step and step-down; gray color)
and the Westfall--Young permutation procedure (single-step and
step-down; dark color).}
\label{figtoeplitz}
\end{figure}

We use a two-sided Wilcoxon test. The power over 250 simulations of the
single-step and step-down methods of the Bonferroni-correction, the
oracle procedure and the Westfall--Young permutation procedure are shown
in Figure \ref{figtoeplitz} for level $\alpha=0.05$. The step-down
version of the Bonferroni correction is the Bonferroni--Holm procedure
\cite{holm79simple} and the step-down version of the Westfall--Young
procedure is given in Westfall and Young \cite{westfall93resampling}.
The oracle
threshold (both single-step and step-down) is approximated on a
separate set of 1000 simulations, and the Westfall--Young method is
using 1000 permutations for each simulation.


The following main results emerge: the Westfall--Young method is very
close in power to the oracle procedure for all values of $m$ in the
block model, giving support to the asymptotic results of Theorem \ref
{thoptimal}. Moreover, the Westfall--Young and the oracle procedures
are also very similar in the Toeplitz model, indicating that
Theorem \ref{thoptimal} may be generalized beyond block-independence
models. The power gains of the Westfall--Young procedure, compared to
Bonferroni--Holm, are substantial; between 20\% and 250\% for the
considered scenarios, where the largest gains are achieved in settings
with a large number of hypotheses and high correlations. Finally, the
difference between step-down and single-step methods is very small in
these sparse high-dimensional settings for all three multiple testing methods.

It might be unexpected that the power of the Westfall--Young is slightly
larger than the oracle procedure for two settings with very high
correlations. This is due to the finite number of simulations when
approximating the oracle threshold, a finite number of permutations in the
Westfall--Young procedure\vadjust{\goodbreak} and a finite number of simulation runs.
The family-wise error rate is between 0.03 and 0.04 for both oracle and
Westfall--Young procedures and below 0.02 for the Bonferroni correction
in the Toeplitz model. The nominal level of $\alpha=0.05$ is exceeded
sometimes in the block model (again for the reason that we use only a
finite number of simulations), where both the oracle and
Westfall--Young procedures attain a family-wise error rate between 0.04
and 0.07 in all settings.


The computational cost of the Westfall--Young procedures scales
approximately linearly with the number $m$ of hypotheses. When using
1000 permutations, computing the Westfall--Young threshold takes about
1.4 seconds per hypotheses on a 3 GHz CPU. For the largest setting of
$m=10\mbox{,}000$ hypotheses, the threshold could thus be computed in just over
2 minutes. It seems as if this computational cost is acceptable, even
for very large-scale testing problems.

\section{Discussion}
We considered asymptotic optimality of large-scale multiple testing
under dependence within a nonparametric framework.
We showed that, under certain assumptions, the Westfall--Young
permutation method is asymptotically optimal in the following sense:
with probability
converging to 1 as the number of tests increases, the
Westfall--Young critical value for multiple testing at nominal level
$\alpha$ is greater
than or equal to the unknown
oracle threshold at level $\alpha-\delta$ for any $\delta>0$.
This implies that the actual level of the Westfall--Young procedure
converges to the
effective oracle level $\alpha_-$.
To investigate the possible impact of discrete $p$-values, we studied a
specific example and provided sufficient conditions that ensure that
$\alpha_-$ converges to $\alpha$.\looseness=-1

We gave several examples that satisfy subset-pivotality and our assumptions
(B1)--(B3) [while assumptions (A1)--(A3) are about the unknown
data-generating distribution].
Most of these examples involve
rank-based or permutation tests.
These tests are
appropriate for very
high-dimensional testing problems. If the number of tests is in the
thousands or even millions, extreme tail probabilities are required
to claim significance, and these tail probabilities
are more
trustworthy under a nonparametric than a parametric test.

If the hypotheses are strongly
dependent, the gain in power
of the Westfall--Young method compared to a simple Bonferroni correction
can be very substantial.
This is a well known, empirical fact, and
we have established here that
this improvement is also optimal
in the asymptotic framework we considered.

Our theoretical results could be expanded to include step-down
procedures like
Bonferroni--Holm \cite{holm79simple}
and the step-down Westfall--Young method \cite{westfall93resampling,GeEtAl03}.
The distinction between single-step and step-down procedures is very
marginal though
in our sparse
high-dimensional framework, as reported in Section \ref{secempirical},
since the number of rejected hypotheses is orders of magnitudes
smaller than the total number of hypotheses.

\section{Proofs}\label{secproofs}
After introducing some additional notation in Section \ref
{secnotationproofs}, the proof of Theorem \ref{thoptimal} is in Section
\ref{sectionproof1} and the proof of Theorem \ref{thoptimaldiscrete} is
in Section~\ref{sectionproof2}.

\subsection{Additional notation}\label{secnotationproofs}
Let $p^{(b)}(W)$ be the minimum $p$-value over all true null hypotheses
in the $b$th block:
\[
p^{(b)}(W) = \min_{j\in A_b\cap I(P_m) } p_{j}(W),\qquad
b=1,\ldots,B,\vadjust{\goodbreak}
\]
and let $\pi_b(c)$ denote the probability 
under $P_m$ that $p^{(b)}(W)$ is less than or equal to a constant
$c\in[0,1]$, that is,
\[
\pi_b(c) = P_{m}\bigl(p^{(b)}(W) \le c\bigr),\qquad b=1,\ldots,B.
\]
Throughout, we denote the expected value, the variance and the
covariance under $P_m$ by $E_m$, $\mathrm{Var}_m$ and $\mathrm{Cov}_m$,
respectively.

\subsection{\texorpdfstring{Proof of Theorem \protect\ref{thoptimal}}{Proof of Theorem 1}}
\label{sectionproof1}
Let $\alpha' \in(0,1)$
and $\delta' \in(0,\alpha')$. Let $\delta=\delta'/2$ and $\alpha
=\alpha'-\delta'$. Then writing expression (\ref{eqthoptimal}) in
terms of $\alpha'$ and $\delta'$ is equivalent to
\[
P_m\{ \hat c_{m,n}(\alpha+2\delta) \ge c_{m,n}(\alpha) \} \to
1 
\qquad \mbox{as } m\to\infty.
\]
By definition,
\[
\hat c_{m,n}(\alpha+2\delta) = \max\Bigl\{ s \in S_n \dvtx P^*\Bigl(\min
_{j\in
\{1,\ldots,m\}} p_{j}(W) \le s\Bigr) \le\alpha+2\delta\Bigr\}.
\]
We thus have to show that
%
\begin{equation}
P_m \Bigl\{ P^*\Bigl(\min_{j\in\{1,\ldots,m\}} p_{j}(W) \le
c_{m,n}(\alpha)\Bigr) \le\alpha+2\delta\Bigr\} \to1
\end{equation}
as $m\to\infty$.

First, we show in Lemma \ref{lemmachatI} that there exists
an $M<\infty$ such that
\[
P^*\Bigl(\min_{j\in\{1,\ldots,m\}} p_{j}(W) \le c_{m,n}(\alpha)\Bigr)
\le P^*\Bigl(\min_{j\in I(P_m)} p_{j}(W) \le c_{m,n}(\alpha)\Bigr)
+\delta
\]
for all $m>M$ and for all $W$.
This result is
mainly due to the sparsity assumption~(A2). Second, we show in Lemma
\ref{lemmaapproxpermutationallsimple} that
%
\begin{equation}
P_m \Bigl\{ P^*\Bigl(\min_{j\in I(P_m)} p_{j}(W) \le c_{m,n}(\alpha)
\Bigr) \le\alpha+\delta\Bigr\} \to1 \qquad\mbox{for } m\to\infty.
\end{equation}
Theorem \ref{thoptimal} follows by combining these two results.
%
\begin{lemma}\label{lemmachatI}
Let $\alpha\in(0,1)$, $\delta\in(0,\alpha)$, and assume
\textup{(A1), (A2), (B2)}\break and~\textup{(B3)}.
Then there exists an $M<\infty$ such that
\[
P^*\Bigl(\min_{j\in\{1,\ldots,m\}} p_{j}(W) \le c_{m,n}(\alpha) \Bigr)
\le P^*\Bigl(\min_{j\in I(P_m)} p_{j}(W) \le c_{m,n}(\alpha) \Bigr)
+\delta
\]
for all $m>M$ and for all $W$.
\end{lemma}
\begin{pf}
Note that $c_{m,n}(\alpha)\in S_n$ by definition.
Using the union bound, we 
have, for all $s\in S_n$ and all $W$,
%
\begin{eqnarray}
P^*\Bigl(\min_{j\in\{1,\ldots,m\}} p_{j}(W) \le s\Bigr) &\le& P^*\Bigl(\min_{j\in
I(P_m)} p_{j}(W) \le s\Bigr)\nonumber\\[-8pt]\\[-8pt]
&&{} + \sum_{j \in I'(P_m)} P^* \bigl(p_j(W) \le s\bigr).\nonumber
\end{eqnarray}
Hence, we only need to show
that there exists an $M<\infty$ such that
%
\begin{equation}\label{eqtoshownew}
\sum_{j \in I'(P_m)} P^* \bigl(p_j(W) \le c_{m,n}(\alpha)\bigr) \le\delta
\end{equation}
for all $m>M$ and all $W$. By assumption (B2)
with constant $\con$,
%
\begin{eqnarray} \label{eqlastline}
\sum_{j \in I'(P_m)} P^* \bigl(p_j(W) \le c_{m,n}(\alpha)\bigr)
&\le&|I'(P_m)| r c_{m,n}(\alpha) \nonumber\\[-8pt]\\[-8pt]
&=& r \frac{|I'(P_m)|}{B} B
c_{m,n}(\alpha).\nonumber
\end{eqnarray}
Since
$|I'(P_m)|/B \to0$ as $m\to\infty$ by assumption (A2), and $B
c_{m,n}(\alpha)$ is bounded above by $-\log(1-\alpha)$
under assumptions (A1) and (B3) (see Lem\-ma~\ref{lemmasumpib}), we
can choose
a $M< \infty$ such that
the right-hand side of
(\ref{eqlastline}) is bounded above
by $\delta$ for all $m>M$.
This proves the claim in (\ref{eqtoshownew}) and completes the proof.
\end{pf}
%
\begin{lemma}\label{lemmaapproxpermutationallsimple}
Let
$\alpha>0$ and $\delta>0$
and assume
\textup{(A1), (A3)} and \textup{(B1)--(B3)}.
Then
\[
P_m \Bigl\{ P^*\Bigl(\min_{j\in I(P_m)} p_{j}(W) \le c_{m,n}(\alpha)
\Bigr) \le\alpha+\delta\Bigr\} \to1
\qquad\mbox{as } m\to\infty.
\]
\end{lemma}
\begin{pf}
Let $\varepsilon>0$.
The statement in the lemma is equivalent to showing that there exists
an $M<\infty$ such that
%
\begin{equation}\label{eqts2}
P_{m} \Bigl\{ P^*\Bigl( \min_{j\in I(P_m)} p_j(W) >c_{m,n}(\alpha)
\Bigr) < 1- \alpha- \delta\Bigr\} <\varepsilon
\end{equation}
for all $m>M$.
By definition,
%
\begin{eqnarray} \label{eqPstareq}\quad
P^*\Bigl( \min_{j\in I(P_m)} p_j(W) > c_{m,n}(\alpha) \Bigr) & = & \frac
{1}{|\mathcal{G}|} \sum_{g\in\mathcal G} 1\Bigl\{ \min_{j\in I(P_m)}
p_j(gW) > c_{m,n}(\alpha) \Bigr\} \nonumber\\[-8pt]\\[-8pt]
& = & \frac{1}{|\mathcal{G}|} \sum_{g \in\mathcal G} R(g,W),\nonumber
\end{eqnarray}
where
\[
R(g,W) := 1\Bigl\{ \min_{j \in I(P_m)} p_j(gW) >c_{m,n}(\alpha) \Bigr\}.
\]
(We
suppress the dependence on $m,n,P_m$ and $\alpha$ for notational simplicity.)

Let $G$ be a random permutation, chosen uniformly in $\mathcal G$, and
let $1$ denote the identity permutation. Then,
by assumption (B1),
it follows that
\[
E_{m} \biggl( \frac{1}{|\mathcal{G}|} \sum_{g \in\mathcal G} R(g,W)
\biggr) = E_{m,G} R(G,W) = E_{m} R(1,W).
\]
By definition of $c_{m,n}(\alpha)$ [see (\ref{eqoracle})],
\[
E_m R(1,W) = P_m \Bigl( \min_{j\in I_{P_m}} p_j(W) > c_{m,n}(\alpha)\Bigr) \ge
1-\alpha.
\]
Hence, the desired result (\ref{eqts2}) follows from a Markov inequality
as soon as one can show that the variance of (\ref{eqPstareq})
vanishes 
as $m\to\infty$, that is, if
%
\begin{equation}\label{eqPstarVar}
\mathrm{Var}_{m} \biggl( \frac{1}{|\mathcal{G}|} \sum_{g \in\mathcal G}
R(g,W) \biggr)\,{=}\,\frac{1}{(|\mathcal{G}|)^2}\!\sum_{g,g' \in\mathcal
G}\!\mathrm{Cov}_{m}(R(g,W), R(g',W) ) = o(1)\hspace*{-35pt}
\end{equation}
as $m\to\infty$.

Let $G,G'$ be two
random permutations, drawn independently and uniformly
from $\mathcal{G}$. Then
\[
\mathrm{Cov}_{m,G,G'} ( R(G,W), R(G',W) ) = \frac{1}{(|\mathcal
{G}|)^2} \sum_{g,g' \in\mathcal G} \operatorname{Cov}_{m}(R(g,W), R(g',W)
).
\]
Hence, in order to show (\ref{eqPstarVar}), we only need to show that
\[
\mathrm{Cov}_{m,G,G'} ( R(G,W), R(G',W) ) =o(1) \qquad\mbox{for
} m\to\infty.
\]
Define
%
\begin{equation}
R_b(g,W) := 1\bigl\{ p^{(b)}(g W) > c_{m,n}(\alpha)\bigr\},
\end{equation}
so that
$R(g,W):= \prod_{b=1}^B R_b(g,W )$.
We then need to
prove that, as $m\to\infty$,
%
\begin{equation}\label{eqts3}
E_{m,G,G'} \Biggl( \prod_{b=1}^B R_b(G,W) R_b(G',W) \Biggr)\,{-}\,\Biggl(
E_{m,G}\Biggl( \prod_{b=1}^B R_b(G,W)\Biggr) \Biggr)^2\,{=}\,o(1) .\hspace*{-35pt}
\end{equation}
Using assumption (A1), the
left-hand side in (\ref{eqts3}) can be written as
\[
\prod_{b=1}^B E_{m,G,G'} \{ R_b(G,W) R_b(G',W) \} - \prod_{b=1}^B
[ E_{m,G} \{R_b(G,W) \} ]^2 .
\]
Note that
$E_{m,G,G'} \{ R_b(G,W)R_b(G',W) \}$
and
$[ E_{m,G}\{ R_b(G,W)\}]^2$ are boun\-ded between $0$ and $1$. For
sequences of numbers $a_1,\ldots,a_B$ and $b_1,\ldots,b_B$ that are
boun\-ded between $0$ and $1$,
the following inequality holds:
\[
\Biggl| \prod_{j=1}^B a_j - \prod_{j=1}^B b_j \Biggr| = \Biggl| \sum
_{j=1}^B \biggl\{ (a_j-b_j)\biggl(\prod_{k<j} b_k \biggr)\biggl(\prod
_{k>j}a_{k}\biggr) \biggr\} \Biggr| \le\sum_{j=1}^B |a_j-b_j
| .
\]
Hence, in order to show (\ref{eqts3})
it is sufficient to show that
%
\begin{eqnarray}\label{eqts4}
&&{\max_{b=1,\ldots,B}}
| E_{m,G,G'} \{ R_b(G,W) R_{b}(G',W)\} - [ E_{m,G} \{
R_b(G,W)\}]^2
| \nonumber\\[-8pt]\\[-8pt]
&&\qquad=
o(B^{-1})\nonumber
\end{eqnarray}
as $m\to\infty$.

Conditional on $W$,
\[
R_b(G,W), R_b(G',W) \stackrel{\mathrm{i.i.d.}}{\sim}
\operatorname{Bernoulli}(\mu_b(W)),\vadjust{\goodbreak}
\]
where $\mu_b(W)$ is the proportion of all
permutations 
$g \in\mathcal G$ for which $R_b(g,W )=1$ or, equivalently,
%
\begin{equation}\label{eqdefmub}
\mu_b(W) = P_{m,G}\bigl\{p^{(b)}(GW) > c_{m,n}(\alpha)|W \bigr\}
= P^*\bigl\{ p^{(b)}(W) > c_{m,n}(\alpha)\bigr\}.\hspace*{-20pt}
\end{equation}
Thus, the
random proportion $\mu_b(W)$ is a function of $W$. Denote its
distribution by~$F_b$.
Using Lemma \ref{lemmasupportF}, the support of $F_b$ is
contained in the interval
$[1 - \log\{1/(1-\alpha)\} \alpha\con^2 m_B B^{-1} , 1]$ under
assumptions (A1), (B1) and (B2). Hence, using
Lemma \ref{lemmavar}, it follows that
\begin{eqnarray*}
0 & \le & E_{m,G,G'}\{ R_b(G,W) R_b(G',W)\} - [ E_{m,G}\{ R_b(G,W)\}
]^2 \\
& \le &\bigl( \log\{1/(1-\alpha)\}\alpha\con^2 m_B B^{-1}\bigr)^2.
\end{eqnarray*}
Since $m_B =o(\sqrt{B})$ under assumption (A3), claim (\ref{eqts4}) follows.
\end{pf}
%
\begin{lemma}\label{lemmasumpib}
Under
assumptions \textup{(A1)} 
and \textup{(B3)}, we have
%
\begin{equation}\label{eqboundsumpib}
Bc_{m,n}(\alpha) \le\sum_{b=1}^B \pi_b(c_{m,n}(\alpha)) \le\log\{
1/(1-\alpha)\}.
\end{equation}
\end{lemma}
\begin{pf}
Let $b\in\{1,\ldots,B\}$ and $j_b \in I(P_m)\cap A_b$. Then
%
\begin{equation}\label{eqpibc}
\pi_b\{ c_{m,n}(\alpha)\} \ge P_m\bigl(p_{j_b}(W) \le c_{m,n}(\alpha
)\bigr) = c_{m,n}(\alpha),
\end{equation}
where the inequality follows from the definition of $\pi_b(\cdot)$, and
the equality follows from assumption (B3) and the fact that
$c_{m,n}(\alpha)\in S_n$. Summing (\ref{eqpibc}) over $b=1,\ldots,B$
yields the first inequality of (\ref{eqboundsumpib}).

To prove the second inequality of (\ref{eqboundsumpib}), 
note that assumption (A1) and the definition of $c_{m,n}(\alpha)$
imply that
%
\begin{equation}\label{eqconstraint}
1-\prod_{b=1}^B [ 1-\pi_b\{ c_{m,n}(\alpha) \} ]\le\alpha.
\end{equation}
The maximum of
$\sum_{b=1}^B \pi_b\{ c_{m,n}(\alpha)\}$ under constraint (\ref
{eqconstraint})
is obtained when
\[
\pi_1\{c_{m,n}(\alpha)\}=\cdots= \pi_B\{c_{m,n}(\alpha)\}.
\]
This implies $\pi_b\{c_{m,n}(\alpha)\} \le1-(1-\alpha)^{1/B}$ for all
$b=1,\ldots,B$, so that
\[
\sum_{b=1}^B \pi_b\{c_{m,n}(\alpha)\} \le B-B(1-\alpha)^{1/B},
\]
and this is bounded above by $-\log(1-\alpha)$ for all 
values of $B$.
\end{pf}
%
\begin{lemma} \label{lemmasupportF}
Assume \textup{(A1), (B1)} and \textup{(B2)}.
Let
$F_b$ be the distribution of~$\mu_b(W)$, where $\mu_b(W)$ is defined in
(\ref{eqdefmub}). Then
\[
\operatorname{support}(F_b) \subseteq[1 - \log\{1/(1-\alpha)\} \alpha\con^2
m_B B^{-1} , 1] .
\]
\end{lemma}
\begin{pf}
Using assumption (B2) with constant $\con$ and the union bound, it
holds that
\[
1-\mu_b(W) =
P^*\bigl\{ p^{(b)}(W) \le c_{m,n}(\alpha) \bigr\}
\le r |A_b| c_{m,n}(\alpha) .
\]
Since $m_B=\max_{b=1,\ldots,B} |A_b|$,
the support of $F_b$ is thus in the 
interval $[1 - m_B \con c_{m,n}(\alpha) , 1]$.

Hence, the proof is complete if we show that
%
\begin{equation} \label{equpperboundc}
c_{m,n}(\alpha) \le-\log(1-\alpha) \alpha\con B^{-1}.
\end{equation}
To see that (\ref{equpperboundc}) holds, we first show that
%
\begin{equation}\label{equpperboundc2}
1-\alpha\le P_{m}\Bigl\{\min_{j\in I(P_m)} p_j(W) > c_{m,n}(\alpha)\Bigr\}
\le\bigl(1-c_{m,n}(\alpha) /\con\bigr)^B.
\end{equation}
The first inequality in (\ref{equpperboundc2}) follows directly
from the definition of $c_{m,n}(\alpha)$; see (\ref{eqoracle}). To
prove the second inequality, note that assumption (A1) implies that
%
\begin{equation}\label{eqindep}
P_{m}\Bigl\{\min_{j\in I(P_m)} p_j(W) > c_{m,n}(\alpha)\Bigr\} = \prod
_{b=1}^B P_{m}\bigl\{ p^{(b)} (W) > c_{m,n}(\alpha)\bigr\}.
\end{equation}
By assumption (B1) and
the law of iterated expectations,
%
\begin{eqnarray} \label{eqpmpbW}
P_{m}\bigl\{ p^{(b)}(W) > c_{m,n}(\alpha)\bigr\} &= & P_{m,G}\bigl\{ p^{(b)}(GW) >
c_{m,n}(\alpha)\bigr\} \nonumber\\[-8pt]\\[-8pt]
&=& E_{m} \bigl\{ P_{m,G}\bigl\{ p^{(b)} (GW) > c_{m,n}(\alpha)|W\bigr\} \bigr\}.\nonumber
\end{eqnarray}
By assumption (B2), the conditional probability within each block
satisfies
%
\begin{eqnarray} \label{eqpmG}
P_{m,G}\bigl\{ p^{(b)}(GW) > c_{m,n}(\alpha)|W\bigr\} & = & P^*\bigl\{ p^{(b)}(W) >
c_{m,n}(\alpha)\bigr\} \nonumber\\
& \le& 1- P^*\{ p_{j_b}(W) \le c_{m,n}(\alpha) \} \\
&\le& 1-c_{m,n}(\alpha) / \con,\nonumber
\end{eqnarray}
where $j_b \in I(P_m)\cap A_b$.
Since the right-hand side of (\ref{eqpmG}) does not depend on~$W$, the
same bound
holds for (\ref{eqpmpbW}), where we also take the expectation over~$W$.
Using this result in (\ref{eqindep}),
the second inequality in (\ref{equpperboundc2}) follows.
Finally,~(\ref{equpperboundc2}) implies
\[
c_{m,n}(\alpha) \le r \{ 1- (1-\alpha)^{1/B} \} .
\]
Since $B(1-(1-\alpha)^{1/B}) \le-\log(1-\alpha)$ for all values of
$B$, it follows that $1-(1-\alpha)^{1/B} \le-\log(1-\alpha) B^{-1}$.
This proves (\ref{equpperboundc}) and completes the proof.~%
\end{pf}
%
\begin{lemma}\label{lemmavar}
Let $U$ be a real-valued random variable with
support $[a,b] \subset[0,1]$. Suppose that the distribution of
the
two random variables $X_1$ and $X_2$, conditional on $U=u$, is given by
\[
X_1,X_2\stackrel{\mathrm{i.i.d.}}{\sim} \operatorname{Bernoulli} (u).
\]
Then $0 \le E(X_1X_2) - E(X_1)E(X_2) \le(b-a)^2$.\vadjust{\goodbreak}
\end{lemma}
\begin{pf}
By the assumption that $X_1$ and $X_2$ are Bernoulli conditional on
$U$, it follows that $E(X_1|U) = E(X_2|U) = U$. Combining this with the
law of iterated expectation and the fact that $X_1$ and $X_2$ are
conditionally independent given $U$, we obtain
\[
E(X_1X_2) = E_U\{ E(X_1X_2|U)\} = E_U \{E(X_1|U) E(X_2|U) \} = E(U^2).
\]
Moreover, we have $E(X_1)=E_U \{E(X_1|U)\} = E(U)$ and similarly
$E(X_2)=E(U)$. Hence,
\[
E(X_1X_2) - E(X_1)E(X_2) = E(U^2) - \{E(U)\}^2 = \operatorname{Var}(U).
\]
Finally, $0 \le\operatorname{Var}(U) \le(b-a)^2$ by the assumption
on the
support of $U$.
\end{pf}

\subsection{\texorpdfstring{Proof of Theorem \protect\ref{thoptimaldiscrete}}{Proof of Theorem 2}}
\label{sectionproof2}
First, note that (W) implies (B1)--(B3).
Using the union bound and assumption (B3), it holds for any $s\in S_n$
that $ms$ is an upper bound for
$ P_m(\min_{j\in I(P_m)} p_j(W) \le s)$.
Hence,
%
\begin{eqnarray} \label{eqcmnge}
c_{m,n}(\alpha) & = & \max\Bigl\{ s \in S_n \dvtx P_m\Bigl( \min_{j\in I(P_m)}
p_{j}(W)\le s\Bigr) \le\alpha\Bigr\} \nonumber\\[-8pt]\\[-8pt]
& \ge & \max\{ s \in S_n \dvtx m s \le\alpha\}.\nonumber
\end{eqnarray}
This implies that the
oracle critical value is larger than zero if the set $\{s \in S_n \dvtx
ms \le\alpha\}$ is nonempty, which is the case if $\min(S_n) \le
\alpha/m$. The smallest possible two-sided Wilcoxon $p$-value is $\min
(S_n) = 2 \frac{(n/2)!(n/2)!}{n!} \le2^{-n/2+1}$.
Hence, it is sufficient to require that $2^{-n/2+1} \le\alpha/m$, or
equivalently, that $n \ge2 \log_2(m/\alpha) + 2$.

Note that (A3$'$) implies (A3). Hence, under assumptions (W), (A1), (A2)
and (A3$'$), the result in Theorem \ref{thoptimal} applies.

Let $\alpha\in(0,1)$.
We will now show that under assumptions (W), (A1) and~(A3$'$),
\[
\alpha_- \to\alpha
\]
as $m,n\to\infty$ such that $n/\log(m) \to\infty$, where $\alpha_-$
was defined in (\ref{eqalpha-}).
Define $c^+_{m,n}(\alpha) := \min\{s\in S_n\dvtx s>c_{m,n}(\alpha)\}$.
Using the definition of $\alpha_-$ and assumption~(A1), we have
%
\begin{eqnarray} \label{eqstart}
\alpha_- & = & P_m\Bigl( \min_{j\in I(P_m)} p_j(W) \le c_{m,n}(\alpha)
\Bigr)\nonumber\\
& = & 1-\prod_{b=1}^B [1-\pi_b\{c_{m,n}(\alpha)\}] \\
& = & 1-\prod_{b=1}^B [1-\pi_b\{c^+_{m,n}(\alpha)\} + \pi_b\{
c^+_{m,n}(\alpha)\} - \pi_b\{c_{m,n}(\alpha)\} ].\nonumber
\end{eqnarray}

Define the function $g_{m,n}\dvtx \prod_{b=1}^B [0,\pi_b\{
c_{m,n}^+(\alpha
)\}] \to\R$ by
\begin{eqnarray*}
g_{m,n}(u) &:=& g_{m,n}(u_1,\ldots,u_B) \\
&:=& 1-\prod_{b=1}^B [1- \pi_b\{
c^+_{m,n}(\alpha)\} + u_b],
\end{eqnarray*}
so that the right-hand side of (\ref{eqstart}) equals $g_{m,n}(w)$,
where $w_b:= \pi_b\{c^+_{m,n}(\alpha)\}-\pi_b\{c_{m,n}(\alpha)\}$ for
$b=1,\ldots,B$. A first-order Taylor expansion of $g_{m,n}(w)$ around
$(0,\ldots,0)$ yields
%
\begin{equation}\label{eqtaylor}
\alpha_- = g_{m,n}(w) = g_{m,n}(0) + \sum_{b=1}^B w_b \,\frac
{\partial g_{m,n}(u)}{\partial u_b}\bigg|_{u=0} + R,
\end{equation}
where $R = o(\sum_{b=1}^B w_b)$.
For all $b=1,\ldots,B$, we have
\begin{eqnarray*}
\frac{\partial g_{m,n}(u)}{\partial u_b}\bigg|_{u=0} & = & - \prod
_{j=1,j\neq b}^B[ 1-\pi_j\{c^+_{m,n}(\alpha)\} ]\\
& = &
- \frac{1-g_{m,n}(0)}{1-\pi_b\{c^+_{m,n}(\alpha)\}} \ge- \frac
{1-g_{m,n}(0)}{1-m_B c^+_{m,n}(\alpha)},
\end{eqnarray*}
where the inequality follows from $\pi_b\{c^+_{m,n}(\alpha)\} \le m_B
c^+_{m,n}(\alpha)$ for $b=1,\ldots,B$, by the union bound and
assumption (B3).
Plugging this into (\ref{eqtaylor}) yields
%
\begin{eqnarray} \label{eqtaylor2}
\alpha_- & \ge & g_{m,n}(0) - \frac{1-g_{m,n}(0)}{1-m_B
c^+_{m,n}(\alpha
)}\sum_{b=1}^B w_b + R \nonumber\\[-8pt]\\[-8pt]
& = & g_{m,n}(0) \biggl( 1 + \frac{\sum_{b=1}^B w_b}{1-m_B
c^+_{m,n}(\alpha)}\biggr) - \frac{\sum_{b=1}^B w_b}{1-m_B
c^+_{m,n}(\alpha)} + R.\nonumber
\end{eqnarray}
The definition of $c^+_{m,n}(\alpha)$ implies that $g_{m,n}(0)>\alpha$
for all $m$ and $n$. Hence, if
%
\begin{equation}\label{eqtwothingstoshow}
\sum_{b=1}^B w_b \to0 \quad\mbox{and}\quad m_Bc^+_{m,n}(\alpha)
\to0
\end{equation}
as $m,n\to\infty$ such that $n/\log(m)\to\infty$, then the right-hand
side of (\ref{eqtaylor2}) converges to $\alpha$ and the proof is complete.

We first consider $\sum_{b=1}^B w_b$. By definition, there is no value
$s'\in S_n$ such that $c_{m,n}(\alpha) < s' < c^+_{m,n}(\alpha)$. Hence,
\begin{eqnarray*}
w_b & = & P_m\Bigl\{ \min_{j \in A_b \cap I(P_m)}p_j(W) = c^+_{m,n}(\alpha
) \Bigr\} \\
& \le & m_B \max_{j\in A_b\cap I(P_m)} P_m\{ p_j(W)=c^+_{m,n}(\alpha)\}
\\
& = & m_B\{c^+_{m,n}(\alpha) - c_{m,n}(\alpha)\},
\end{eqnarray*}
where the inequality follows from the union bound, and the last
equality is due to assumption (B3). This implies
\[
\sum_{b=1}^B w_b \le B m_B \{c^+_{m,n}(\alpha) - c_{m,n}(\alpha)\} = B
c_{m,n}(\alpha) m_B \biggl( \frac{c^+_{m,n}(\alpha)}{c_{m,n}(\alpha)} -
1\biggr).
\]
Similarly, we have
\[
m_B c^+_{m,n}(\alpha) = B c_{m,n}(\alpha) \frac{m_B}{B} \frac
{c^+_{m,n}(\alpha)}{c_{m,n}(\alpha)}.
\]
Note that $B c_{m,n}(\alpha) \le\log\{1/(1-\alpha)\}$ by Lemma \ref
{lemmasumpib} and $m_B=O(1)$ (and hence $B\to\infty$) by assumption
(A3$'$). Hence, in order to prove (\ref{eqtwothingstoshow}), it
suffices to show that
%
\begin{equation} \label{eqratioConv}
c_{m,n}^+(\alpha)/c_{m,n}(\alpha) \to1 \qquad\mbox{as } m,n \to\infty,
n/\log(m) \to\infty.
\end{equation}
Let the ordered $p$-values in $S_n$, based on a two-sided Wilcoxon test
with equal sample sizes $n/2$ in both classes, be denoted by $s_0 < s_1
< \cdots<s_{r_n}$, where $r_n = \lfloor n^2/8 + 1 \rfloor$. It is well
known that
\[
s_i = 2 \frac{(n/2)!(n/2)!}{n!} \sum_{j=0}^{i} q_{n/2}(j) \qquad\mbox
{for } i=0,\ldots,r_n-1
\]
and $s_{r_n}=1$, where $q_n(j)$ is the number of integer partitions of
$j$ such that neither the number of parts nor the part magnitudes
exceed $n$ \mbox{[and $q_n(0)\,{=}\,1$]}~\cite{Wilcoxon45}. Let $i_{m,n}$ satisfy
$s_{i_{m,n}} = c_{m,n}(\alpha)$. Then
\[
\frac{c^+_{m,n}(\alpha)}{c_{m,n}(\alpha)} = \frac{\sum
_{j=0}^{i_{m,n}+1} q_{n/2}(j)}{\sum_{j=0}^{i_{m,n}} q_{n/2}(j)}.
\]
This ratio converges to 1 if $i_{m,n}\to\infty$.
Recall that $c_{m,n}(\alpha) \ge\max\{s\in S_n \dvtx s \le\alpha
/m\}$
[see (\ref{eqcmnge})]. Hence,
%
\[
c^+_{m,n}(\alpha) = 2 \frac{(n/2)!(n/2)!}{n!} \sum_{j=0}^{i_{m,n}+1}
q_{n/2}(j) > \alpha/m.
\]
Since $2m\{(n/2)!(n/2)!\}/n! \le m 2^{-n/2} \to0$ as $m,n \to\infty$
such that $n/\break\log(m)\to\infty$, we have
that under these conditions $i_{m,n}\to\infty$ and $c^+_{m,n}(\alpha)
/\break c_{m,n}(\alpha) \to1$. Thus (\ref{eqratioConv}) holds and hence
implies (\ref{eqtwothingstoshow}), which completes the proof.

\section*{Acknowledgments}

We would like to thank two referees for constructive comments.


%

\printaddresses

\end{document}